\documentclass[9pt,twocolumn,twoside]{osajnl}
\usepackage{svg}
\usepackage{mathtools}
\usepackage{amsthm}
\usepackage{amssymb}
\journal{ao} 

\setboolean{shortarticle}{false}

\title{Solitary solution method for incompressible Navier-Stokes PDE in $\mathbb{R}^n$}

\author[1]{Johannes Lawen}

\affil[1]{Institute of Process Systems Engineering, Hamburg University of Technology, Am Schwarzenberg-Campus 4 (C), Hamburg, 21073}

\affil[*]{Corresponding author: jl@environment.report}

\begin{abstract}
The method exploits the contraction of space to systematically obtain compact solitary solutions. The latter are provided for the incompressible Euler and Navier-Stokes PDE. The nonlinear response of momentum advection is moved into a term for contracting space. Then the linear continuity PDE is solved by means of arbitrarily selected closure functions. The contracting space is then split into two variables. The compactness of some solutions is enhanced by numerically integrating the contracting domain while retaining a solution for the nonlinear PDE. The validation of numerical schemes is demonstrated for the Euler and Navier-Stokes PDE. As the nonlinear response is isolated in only one spatial dimension, the method permits to validate arbitrary unstructured meshes and domain geometries by introducing the spatial dimension $n+1$. 
\end{abstract}

\setboolean{displaycopyright}{true}

\begin{document}

\maketitle

\section{Introduction}
Nonlinear PDE for fluid systems constrain the space-time geometry of compact analytical solutions which are of valuable validatory utility for numerical solvers \cite{Bristeau2018VariousAS}. Shaping domain geometries such as to attain analytical flow profiles is, hence, plausibly key to obtain a generalized approach: that is, an analytical solution where the nonlinear property of the velocity is mitigated by moving it into a spatially adapted and dynamic domain geometry. Solitary solutions with hydrodynamic constituents exhibiting continuous wave patterns have been found prior \cite{sinha_long_2014, leoni2020traveling, drazin_solitons:_1989}.
\newline
A dynamic fluid domain geometry is in principle not alien to free surface flow, given, for example, tidal ocean dynamics and evolving seabeds that are subject to sediment fluxes. Besides this geophysical example, the general property of a contracting domain or space is also intrinsic to the fabric of relativistic systems as fundamental principle. 
In recent years a cross-utility of such correspondence is increasingly receiving attention \cite{Bredberg2012, Padmanabhan:2010rp, Hubeny_2011, Rangamani_2009}, its utility for nonrelativistic fluid dynamics has been reported \cite{ELING2009496}, and its mutual fertilization is evident in the subject matter of journals such as \emph{Geophysical \& Astrophysical Fluid Dynamics}. 
\newline
The onward sections illustrate how the inconvenience of the nonlinear response of momentum advection can be moved into a term for contracting space. Whereas analytical solutions have been provided before as propositions \cite{Bristeau2018VariousAS}, here a general method is promulgated to obtain such. Also, a composite treatment is demonstrated to deploy analytical solutions for validations of momentum transport by numerically integrating the dynamic domain topography. This averts otherwise inconvenient analytical expressions for the latter. The solution scheme is first demonstrated for one spatial dimension, in section \ref{sec:ndim} for $n$ spatial dimensions, and in two steps: 1. The term
\[\frac{\partial u(x,t)^2}{\partial x}\]
for the velocity $u$ (unit: $m\: s^{-1}$) and the spatial coordinate $x$ (unit: $m$) is denoted with
\[\frac{\partial \left( h(x,t) u(x,t)^2 \right)}{\partial x}\]
absent the skipping of the cross-section $h(x,t)$ (unit: $m^2$) in the infinitesimal balance, as it is kept variable along $x$ (unit: $m$) and time $t$ (unit: $s$). 
In this approach first the linear continuity PDE is solved by means of arbitrarily selected closure functions. 2. The contracting space is split into two variables to prevent the system from becoming overdetermined:

\begin{equation}
    h(x,t) = h_B(x,t) + h_E(x,t)
\end{equation}

E.g. $h_E$ may be the hydrostatic elevation and $h_B$ the topography or, in aquatic terms, the bathymetry.
Henceforth, the Euler PDE is now 
\begin{enumerate}[label=(\roman*), nosep, topsep=8pt]
    \item incompressible,
    \item depicting a fluid with constant density,
    \item and configured for free surface flow to be solved together with the continuity PDE.
\end{enumerate}
The system is still artificial as friction at the bathymetric boundary is not considered. The solutions for $u(x,t)$ and $h(x,t)$ are inserted into the momentum transport PDE where e.g. the hydrostatic pressure gradient is given with $\partial h_E(x,t)/\partial x$ which can be denoted as $\partial \left(h(x,t) - h_B(x,t)\right)/\partial x$.
The next two sections show that it is trivial to select $h(x,t)$ such that it can be solved for $h_B(x,t)$.
But even if it is selected such that the equation cannot be conveniently solved, knowledge of $h_B(x,t)$ is only relevant for the numerical solver. The analytical solution for validation purposes has been determined apriori for $u(x,t)$ and $h(x,t)$. Therefore, even in cases where $u(x,t)$ is selected such that an integration and solution for $h_B(x,t)$ is not convenient or cannot be attained, it can simply be integrated numerically as $h_B(x,t)$ and its approximation error will enter the numerical but not the analytical solution.
Henceforth, this method reliably produces a solution that suffices for validatory purposes:
\begin{enumerate}[label=(\roman*), nosep, topsep=8pt]
    \item The solution of the linear continuity and arbitrarily selected linear closure function is linear and, hence, can be obtained with e.g. Laplace transformation.
    \item The resulting expression for the derivative of the splitting variable $h_B$ can be, if of unknown solution, integrated numerically as it will only constitute a utility in the numerical solution.
\end{enumerate}
An analytical solution for the momentum along with an ODE for the bathymetry is provided in the subsequent Section \ref{sec:2euler}. Its application in the validation of numerical solutions is provided in Section \ref{sec:soliton}. An analytical expression for the bathymetry is given in Section \ref{sec:eulerpde}.

\section{Euler PDE}
\label{sec:2euler}
Provided in this section is an analytical solution for momentum transport and an ODE for the domain $h_B(x,t)$. Figure \ref{fig:contract} below shows distributions for different $c_1$, the numerical and analytical solution for the velocity (right), and the numerically and analytically obtained bathymetry (left). 
\vspace{4pt}
\begin{figure}[htbp]
\centerline{    
\includegraphics[width=1.1\linewidth]{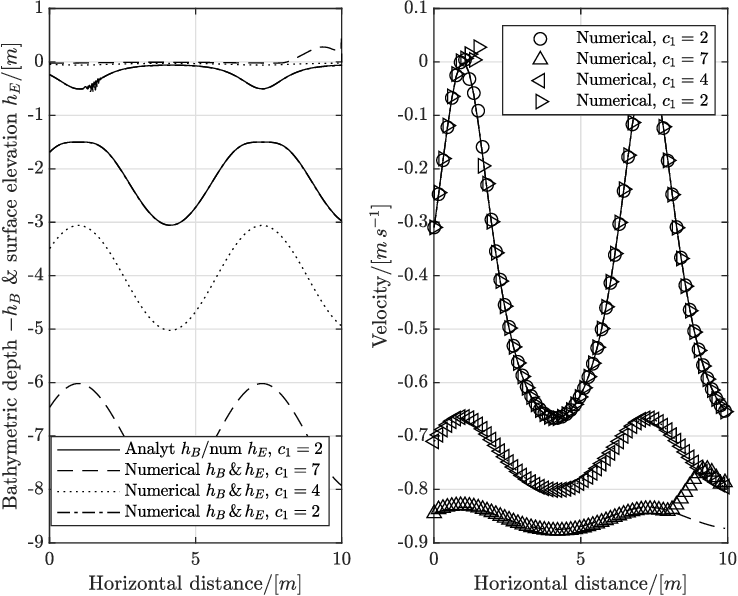}}
  \caption{Analytical and numerical solutions of the Euler momentum  PDE with hydrostatic term, a gravitational acceleration of $1\: m\: s^{-2}$, and the continuity PDE.}
  \label{fig:contract}
\end{figure}
\newline
Input data have been set as per Table \ref{tab:euler-in}. The right side contains the numerical solutions (marked) as per Section \ref{sec:soliton} and the anaytical solutions (lined) for the velocity. All but one solutions are based on the numerical integration of the domain. One solution with analytical bathymetry, as per Section \ref{sec:eulerpde}, has been added for comparison.

\begin{table}[htbp]
\centering
\caption{\bf Grid for Euler PDE example}
\begin{tabular}{ccc}
\hline
Quantity & Increments & Increment size \\
\hline
Distance & $10^{3}$ & $10^{-2}[m]$ \\
Time & $10^{4}$ & $10^{-3}[s]$ \\
$c_1$ & $3$ & $[2\:4\:7]$ \\
\hline
\end{tabular}
  \label{tab:euler-in}
\end{table}

The continuity and Euler PDE with cross-section are:
\[\frac{\partial\left(h(x,t) u(x,t)\right)}{\partial t} =\]
\begin{equation}
    \label{moment}
     -\frac{\partial \left(h(x,t) u(x,t)^2\right)}{\partial x}-h(x,t) g\frac{\partial h_E(x,t)}{\partial x} 
\end{equation}

and the continuity
\begin{equation}
    \label{continuity}
    \frac{\partial h(x,t)}{\partial t} = -\frac{\partial \left(h(x,t) u(x,t)\right)}{\partial x}    
\end{equation}
and $h(x,t) = h_E(x,t)+h_B(x,t)$ holds for the surface elevation and bathymetric depth respectively. 
\vspace{4pt}
\newline
\textit{Proposition 1: For some $t_0 \in \mathbb{R}_+$, let us consider the functions $h$, $h_E$, $h_B$, $u$ defined for $t \geq t_0$ by selecting $h(x,t) = c_1 + \sin(x+t)$ which can be satisfied with a dynamic bathymetry $h_B$.} 
\vspace{4pt}
\newline
\textit{Proof}. The proposition yields for the continuity:
\begin{equation}
    \frac{\partial \left(c_1+\sin(x+t)\right)}{\partial t} = -\frac{\partial \left(\left(c_1+\sin(x+t)\right) u(x,t) \right)}{\partial x}
\end{equation}

Here it can be deduced for the velocity $u(x,t)$ that if
\begin{equation}
    \label{eq:contsol}
    u(x,t) = \left(c_1+\sin(x+t)\right)^{-1} - 1
\end{equation}
then LHS and RHS match and the PDE must be fulfilled.
$h(x,t)$ and $u(x,t)$ can be inserted into the Euler PDE:
\[-\frac{\partial \left(\left(c_1+\sin(x+t)\right) \left(\left(c_1+\sin(x+t)\right)^{-1} - 1\right)^2 \right)}{\partial x}\]
\[-\left(c_1+\sin(x+t)\right)g\frac{\partial h_E(x,t)}{\partial x}\]
\begin{equation}
     =\frac{\partial \left(\left(c_1+\sin(x+t)\right)\left(\left(c_1+\sin(x+t)\right)^{-1} - 1\right) \right)}{\partial t}
\end{equation}
It is taken advantage of the fluid column being split into surface elevation $h_E(x,t)$ and bathymetric fluid depth $h_B(x,t)$
\[-\frac{\partial sin(x+t) }{\partial t} =-\frac{\partial \left(c_1+\sin(x+t)+ \left(c_1+\sin(x+t)\right)^{-1} - 2 \right)}{\partial x}\]
\begin{equation}
     -\left(c_1+\sin(x+t)\right)g\frac{\partial \left(h(x,t)-h_B(x,t)\right)}{\partial x}
\end{equation}
\[-\cos(x+t) =-\left(\cos(x+t)-\frac{cos(x+t)}{\left(c_1+\sin(x+t)\right)^2}\right)\]
\begin{equation}
     -g\left(c_1+\sin(x+t)\right)\left(\cos(x+t)-\frac{\partial h_B(x,t)}{\partial x}\right)
\end{equation}
\begin{equation}
     \label{eq:h_b}
     \frac{\partial h_B(x,t)}{\partial x}=\cos(x+t)-\frac{cos(x+t)}{g\left(c_1+\sin(x+t)\right)^3}
\end{equation}

$h_B(x,t)$ is merely obtained through integration as exercised in Section \ref{sec:eulerpde} on the Euler PDE. 
\qed
\vspace{4pt}
\newline
Even if the complexity of the PDE is migrated into a contracting dimension and then disposed by means of variable splitting, as here for a transient bathyemtry, the numerical integration of the domain generally retains the analytic momentum expression and, hence, the possibility to validate approximations of the same. Henceforth, the latter permits to validate more boundary conditions with an analytical solution.

\section{Validation utility}
\label{sec:soliton}

Before proceeding it may be denoted, that the concept has also been approached in focus of the attained favorable property of constant propagation under the terms solitons and solitary waves \cite{sinha_long_2014, leoni2020traveling}. As per some definitions \cite{drazin_solitons:_1989} the term soliton would require also unimpeded propagation through collisions. The analytical solution for momentum obtained above permits to validate numerical solutions. For the latter the incompressible Euler PDE with continuity and surface elevation will be rearranged into the convective form. For this purpose, first the continuity PDE is denoted in its convective form. The product rule is applied to the continuity PDE to obtain:
\begin{equation}
    \frac{\partial h}{\partial t} = -\left(u(x,t)\frac{\partial h(x,t)}{\partial x} + h(x,t)\frac{\partial u(x,t)}{\partial x} \right)
\end{equation}
Likewise the product rule for derivation is applied to the Euler PDE's time derivative and advective momentum transport term, yielding:
\[-\left( u(x,t)^2\frac{\partial h(x,t)}{\partial x} + 2\:h(x,t)u(x,t)\frac{\partial u(x,t)}{\partial x}\right )\]
\[- h(x,t)\: g\frac{\partial h_E(x,t)}{x}\]
\begin{equation}
     = u(x,t)\frac{\partial h(x,t)}{\partial t} + h(x,t)\frac{\partial u(x,t)}{\partial t}
\end{equation}
Subsequently, the convective form of the continuity PDE is inserted into the Euler PDE which eliminates several terms and recovers the familiar form of the Euler PDE's momentum transport term for stagnant domain geometries or, here, stagnant surface elevations in free surface flow absent the continuity PDE:
\begin{equation}
    \frac{\partial u(x,t)}{\partial t} = -u(x,t)\frac{\partial u(x,t)}{\partial x} - \: g\frac{\partial h_E(x,t)}{x} 
\end{equation} 

This PDE can be subjected to the numerical solver and the result validated with the aforementioned analytical solution. Only the continuity PDE remains absent of an explicit expression and is numerically solved to obtain the evolution of the dynamic domain shape. For convenience a gravitational acceleration of $1 m\: s^{-2}$ is set. The boundary conditions for the system are:
\begin{equation}
    u(0,t) = \left(c_1 + \sin(t)\right)^{-1} - 1    
\end{equation}
\begin{equation}
    u(x_n,t) = \left(c_1 + \sin(x_n+t)\right)^{-1} - 1
\end{equation}
\begin{equation}
    h(0,t) = c_1 + \sin(t)
\end{equation}
\begin{equation}
    h(x_n,t) = c_1 + \sin(x_n+t)
\end{equation}
\begin{equation}
    h_E(0,t) = -\left(c_1 + \sin(t)\right)^{-2}\left(2\,g\right)^{-1}
\end{equation}
\begin{equation}
    h_E(x_n,t) = -\left(c_1 + \sin(x_n+t)\right)^{-2}\left(2\,g\right)^{-1}
\end{equation}
\begin{equation}
    h_B(0,t) = c_1 + \sin(t) - h_E(0,t)
\end{equation}
\begin{equation}
    h_B(x_n,t) = c_1 + \sin(x_n+t) - h_E(x_n,t)
\end{equation}
with the bathymetry being computed numerically as per the above equation for $\partial h_B/\partial x$.
The numerical integration of the artificial bathymetry occurs here in a simple first order form, that is, as summation of products of local derivative and spatial increment. This restricts also the spatial step to lower limits which in turn yields even smaller time steps. Again, this particular expression for the bathymetry is solved analytically in the subsequent Section \ref{sec:eulerpde}.





\section{Analytical contraction}
\label{sec:eulerpde}
The prior sections illustrated that the conjunction of the retention of cross-sections, trigonometric approaches, and variable splitting permit to dispose much of the system's complexity into a fluctuating domain, the flow cross-section or bathymetry. Whereas the above retains an analytical solution for correlation purposes, it does not explicate an expression for the water depth. The explicit expression is given below. Absent the utilization of e.g. the Laplace transformation, simple solutions for $h(x,t)$ and $u(x,t)$ can be found by selecting first $h(x,t)$ in the continuity expression for vertically dynamic domains. That is not surprising as $h(x,t)$ occurs twice and $u(x,t)$ only once in the continuity expression. The solutions are, therefore, simpler if $h(x,t)$ is selected first instead of $u(x,t)$.
\vspace{4pt}
\newline
\textit{Corollary 1: $h_B$ from Proposition 1 can be explicitly expressed.}
\vspace{4pt}
\newline
\textit{Proof.}
Given equations \ref{moment}, \ref{continuity}, a bathymetry compliant to equation \ref{eq:h_b} has been found to permit a trigonometric solution for the velocity. The expression for $h_B$ is structured into summands, permitting the exploitation of the sum rule to obtain the anti-derivative. The anti-derivatives are obtained from equation \ref{eq:h_b} for each term and $h_E(x,t)$ is eliminated with $h=h_E+h_B$, yielding the following solution:

\begin{equation}
    h_B(x,t)=c_1+\sin(x+t)+\left(2\: g\right)^{-1}\left(c_1+\sin(x+t)\right)^{-2}
\end{equation}

\begin{equation}
    h_E(x,t)=-\left(2\: g\right)^{-1}\left(c_1+\sin(x+t)\right)^{-2}
\end{equation}

together with $u(x,t)=\left(c_1+\sin(x+t)\right)^{-1}-1$ and $h(x,t)=c_1+\sin(x+t)$. 
\qed
\vspace{4pt}
\newline

An analytical solution is shown in Figure 1 for $c_1=2$ besides numerical solutions for both, the analytical bathymetry above and a numerical approximation of the bathymetry. In both cases, the solution for the velocity is correlated with a numerical solution, demonstrating its utility in validations. For the latter an explicit upwind approximation has been used for all velocity vector transport terms and an explicit central difference approximation for all other transport terms. Both, the continuity and momentum transport PDE have been brought into the convective form. That is, the conservation form has been detangled into two simpler terms by application of the product rule. This does not ascertain quantity conservation but tends to improve stability.



\section{Navier-Stokes PDE}

For the solution of the incompressible Navier-Stokes equation a term for eddy diffusive dissipation is added to the above Euler PDE. Henceforth, the Navier-Stokes PDE is now 
\begin{enumerate}[label=(\roman*), nosep, topsep=8pt]
    \item incompressible,
    \item depicting a fluid with constant density and eddy viscosity,
    \item and configured for free surface flow to be solved together with the continuity PDE.
\end{enumerate}
Again, solutions are obtained by governing  the fluid column $h(x,t)$ with a simple trigonometric function.
\[\overbrace{\frac{\partial\left(h(x,t) u(x,t)\right)}{\partial t} + \frac{\partial \left(h(x,t) u(x,t)^2\right)}{\partial x}}^{material\: derivative}\]
\begin{equation}
     =\underbrace{k_H\frac{\partial \left(h(x,t)\partial u(x,t)/\partial x\right)}{\partial x}}_{eddy\: viscosity}-\underbrace{h(x,t)\:g\frac{\partial \left(h_E(x,t)\right)}{\partial x}}_{hydrostatic\:pressure}     
\end{equation}
with the horizontal eddy viscosity $k_H$ (unit: $m^2\: s^{-1}$).
\vspace{4pt}
\newline
\textit{Proposition 2: For some $t_0 \in \mathbb{R}_+$, let us consider the functions $h$, $h_E$, $h_B$, $u$ defined for $t \geq t_0$ by selecting $h(x,t) = c_1 + \sin(x+t)$ which can be satisfied with a dynamic bathymetry $h_B$.} 
\vspace{4pt}
\newline
\textit{Proof}.
As for the Euler PDE, the continuity PDE returns equation \ref{eq:contsol}. $h(x,t)$ and $u(x,t)$ can be inserted into the NS PDE configured for incompressible fluid, constant density and eddy viscosity, and free surface flow:
\[\frac{\partial \left(\left(c_1+\sin(x+t)\right)\left(\left(c_1+\sin(x+t)\right)^{-1} - 1\right) \right)}{\partial t}\]
\[=k_H\frac{\partial \left(\left(c_1+\sin(x+t)\right)\partial\left(\left(c_1+\sin(x+t)\right)^{-1} - 1\right)/\partial x\right)}{\partial x}\]
\[-\frac{\partial \left(\left(c_1+\sin(x+t)\right) \left(\left(c_1+\sin(x+t)\right)^{-1} - 1\right)^2 \right)}{\partial x}\]
\begin{equation}
     -\left(\sin(x+t)+c_1\right)g\frac{\partial h_E(x,t)}{\partial x}
\end{equation}

\begin{figure}[htbp]
    \centerline{
    \includegraphics[width=1.1\linewidth]{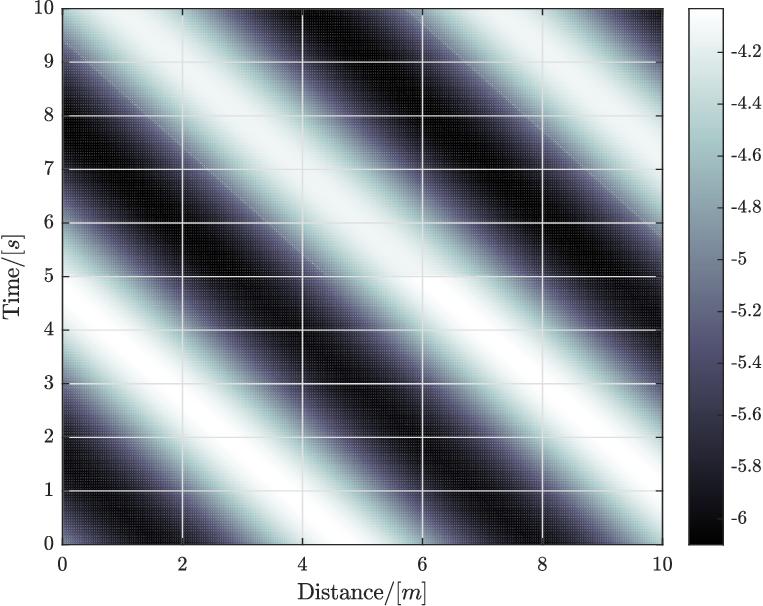}
    }
    \caption{Bathymetry, first $5\,s$ compliant to the solution of the configurations of the Euler PDE and the second $5\,s$ compliant to the NS PDE with $c_1 = 5$ and a large viscosity of 1 $m^2s^{-1}$ for illustrative purposes.}
    \label{fig:NSbathy}
\end{figure}

The dynamic adaptive scale or bathymetry to accomodate the analytic solution is shown in Figure \ref{fig:NSbathy}.
\vspace{4pt}

Again, it is taken advantage of the fluid column being split into surface elevation $h_E(x,t)$ and bathymetric fluid depth $h_B(x,t)$.
\[-\frac{\partial sin(x+t) }{\partial t} =-\frac{\partial \left(c_1+\sin(x+t)+ \left(c_1+\sin(x+t)\right)^{-1} - 2 \right)}{\partial x}\]
\begin{equation}
     -k_H\frac{\partial \frac{\cos(x+t)}{\sin(x+t)+c_1} }{\partial x}-\left(\sin(x+t)+c_1\right)g\frac{\partial h(x,t)-h_B(x,t)}{\partial x}
\end{equation}
\[-\cos(x+t) =-\left(\cos(x+t)-\frac{cos(x+t)}{\left(c_1+\sin(x+t)\right)^2}\right)\]
\[+k_H\frac{\sin(x+t)\left(\sin(x+t)+c_1\right)+\cos^2(x+t)}{\left(\sin(x+t)+c_1\right)^2}\]
\begin{equation}
     -\left(\sin(x+t)+c_1\right)g\left(\cos(x+t)-\frac{\partial h_B(x,t)}{\partial x}\right)
\end{equation}
As the expression above is structured into summands it permits ample exploitation of the sum rule to obtain the anti-derivative.

\[\frac{\partial h_B(x,t)}{\partial x}\]
\[=\cos(x+t)-\frac{k_H}{g}\frac{\sin(x+t)\left(\sin(x+t)+c_1\right)+\cos^2(x+t)}{\left(\sin(x+t)+c_1\right)^3}\]
\begin{equation}
    -\frac{1}{g}\frac{cos(x+t)}{\left(c_1+\sin(x+t)\right)^3}
\end{equation}

Integration yields:
\[h_B(x,t)\]
\[=c_1+\frac{k_H}{g}\left(a(\mathbf{r})+\left(c_1+\sin(x+t)\right)^{-2}\right)\]
\begin{equation}
    +\sin(x+t)+\frac{1}{2g}\left(c_1+\sin(x+t)\right)^{-2}
\end{equation}
with $a(\mathbf{r})$ provided in equation \ref{eq:a} with $\mathbf{r} = [x, t]$ for the one dimensional case.
\[h_E(x,t)=-\frac{1}{g}\left(\frac{\left(c_1+\sin(x+t)\right)^{-2}}{2}\right.\]
\begin{equation}
    \left.+k_H\left(a(\mathbf{r})+\left(c_1+\sin(x+t)\right)^{-2}\right)\right)
\end{equation}

together with $u(x,t)=\left(c_1+\sin(x+t)\right)^{-1}-1$ and $h(x,t)=c_1+\sin(x+t)$. 
\qed
\vspace{4pt}
\newline
For the purpose of the numerical solution the conservation form is replaced with the convective form for the transport terms. 
Applying the product rule on the time derivative and advective transport terms, as well as insertion of the continuity PDE into the NS PDE, as precedingly exercised for the Euler PDE, return:
\[h(x,t)\frac{\partial u(x,t)}{\partial t} = \frac{\partial \left(k_H h(x,t)\partial u(x,t)/\partial x\right)}{\partial x}\]
\begin{equation}
     -h(x,t)u(x,t)\frac{\partial u(x,t)}{\partial x} - g\frac{\partial h_E(x,t)}{x} 
\end{equation} 

Only the continuity PDE remains in the numerical solution reflective of the dynamic domain shape. Only the boundary conditions for the surface elevation differ from those provided for the Euler PDE in the sections above. The boundary conditions for the surface elevation with a constant $k_H$ is for this solution given with:
\[h_E(x,0)=-\frac{1}{g}\left(\frac{\left(c_1+\sin(x)\right)^{-2}}{2}\right.\]
\begin{equation}
    \left.+k_H\left(\left.a(\mathbf{r})\right|_{x=0}+\left(c_1+\sin(x)\right)^{-2}\right)\right)
\end{equation}
\[h_E(x_n,t)=-\frac{1}{g}\left(\frac{\left(c_1+\sin(x_n+t)\right)^{-2}}{2}\right.\]
\begin{equation}
    \left.+k_H\left(\left.a(\mathbf{r})\right|_{x=x_n}+\left(c_1+\sin(x_n+t)\right)^{-2}\right)\right)
\end{equation}
The solutions used for the NS PDE can be arbitrarily set as for the Euler PDE. Just the required seafloor evolution has to change vis-à-vis the case of the Euler PDE to comply with the analytic solution. The difference in the bathymetry evolution is illustrated in Figure \ref{fig:NSbathy}.
The first 5 seconds show the bathymetry without eddy viscosity, that is, the Euler PDE, and the second 5 seconds show a for illustration purposes a large eddy viscosity of 1 $m^{2}/s$, that is, the Navier-Stokes PDE.
\vspace{4pt}
\begin{figure}[htbp]
    \centerline{
    \includegraphics[width=1.1\linewidth]{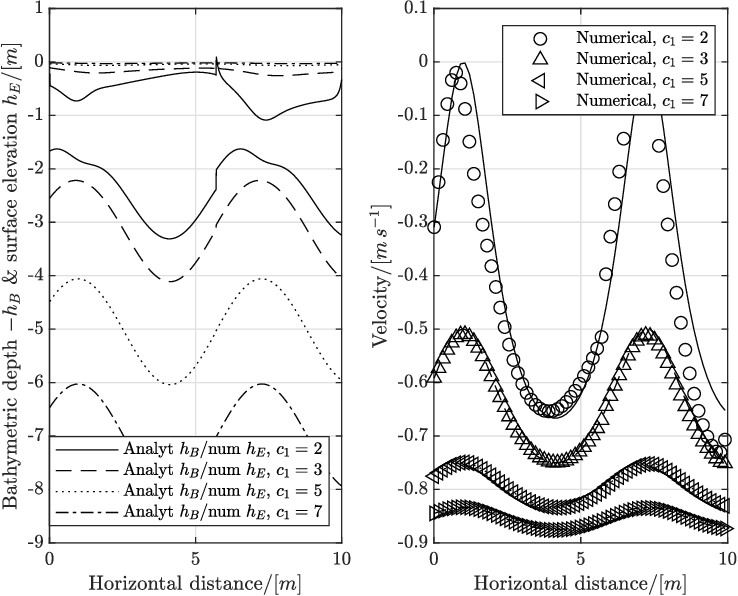}}
    \caption{Analytical and numerical solutions for the velocity of the Navier-Stokes PDE with hydrostatic term, the continuity PDE, and analytical bathymetry with different $c_1$.}
    \label{fig:nsu}
\end{figure}

The simulations for the Navier-Stokes PDE and the presented corresponding analytical solutions were conducted for input data as per Table \ref{tab:ns-in} and a viscosity of $0.3\: m^2\: s^{-1}$.  
\begin{table}[htbp]
\centering
\caption{\bf Grid for Navier-Stokes PDE example}
\begin{tabular}{ccc}
\hline
Quantity & Increments & Increment size \\
\hline
Distance & $10^{3}$ & $10^{-2}[m]$ \\
Time & $10^{5}$ & $10^{-4}[s]$ \\
$c_1$ & $4$ & $[2\:3\:5\:7]$ \\
\hline
\end{tabular}
  \label{tab:ns-in}
\end{table}
\newline
The oscillating distribution shifts in the coordinate system and also changes its amplitude depending on the eddy viscosity. The simulated velocity, shown for $t = 10\: s$ in Figure \ref{fig:nsu}, exhibits a stability dependency to the setting of $c_1$.
\section{n Dimensions}
\label{sec:ndim}
The method can be applied to systems with $n$ spatial dimensions with $U_i$ where $i \in \mathbb{N}$. The NS-PDE, configured for a scalar gravitational field imposed on free hypersurface flow with eddy viscosity, is thus denoted with

\begin{equation}
     \underbrace{\nabla.\left(h\:\mathbf{U }\otimes\mathbf{U }\right)}_{material\:derivative}=\underbrace{\nabla.\left(h\:\nabla.\mathbf{U}\right)}_{eddy\:viscosity} + \underbrace{h\: g \nabla h_E}_{hydrostatic\:pressure} 
\end{equation}

with $\mathbf{U}(\mathbf{r})=[U_1 ... U_{n+1}]^T$, SI units, and the component velocity along time being normalized to one. An individual contracting scale for each spatial dimensions, with $U_i \perp H_i$, is not required for obtaining some solutions. Note, in this notation, henceforth, no steady state system is denoted. All dimensions have a velocity, including time, e.g. $U_t, U_x, U_y, U_z$ with $U_t = 1$. In case of incompressible free surface flow the hydrostatic gradient along a vertical coordinate $z$ would be rendered $\partial H_E/\partial z = 0$. 
The continuity is given with 
\begin{equation}
     \nabla.\left(h\, \mathbf{U }\right)=0 
\end{equation}
The solution approach remains the same:
\begin{enumerate}[label=(\roman*), nosep, topsep=8pt]
\item select a trigonometric solution for $h(\mathbf{r})$ with $\mathbf{r} = [x,y, ... t]$ apriori,
\item select solutions for the velocities that comply with the continuity PDE,
\item split $h(\mathbf{r})$ into $h_B(\mathbf{r})$ and $h_E(\mathbf{r})$, substitute $h_E(\mathbf{r})$ in the momentum PDE, and
\item solve the momentum PDE for $h_B(\mathbf{r})$.
\end{enumerate}
\vspace{4pt}
\textit{Proposition 3: For some $t_0 \in \mathbb{R}_+$, let us consider the functions $h$, $h_E$, $h_B$, $\mathbf{U}$ defined for $t \geq t_0$ by selecting $h(\mathbf{r}) = c_1 + \sin(\sum_ir_i)$ which can be satisfied with a dynamic $h_B$.} 
\vspace{4pt}
\newline
\textit{Proof}.
For the analogon to the prior solution that is, $h(\mathbf{r}) = c_1 + \sin(\sum_ir_i)$, a solution for the continuity PDE can e.g. be symmetrically constructed with $U_i(\mathbf{r}) = \left(\sin(\sum_ir_i)+c_1\right)^{-1} - n^{-1}$ but $U_t = 1$. $h(x,t)$ and $u(x,t)$ can be inserted into the momentum PDE for one component velocity $U_j$:
\[\frac{\partial \left(\left(\sin(\sum_ir_i)+c_1\right)\left(\left(\sin(\sum_ir_i)+c_1\right)^{-1} - n^{-1}\right) \right)}{\partial t}\]
\[=k_H\sum_i\left(\frac{\partial \left(\left(\sin(\sum_ir_i)+c_1\right)\partial\left(\left(\sin(\sum_ir_i)+c_1\right)^{-1} - n^{-1}\right)/\partial r_i\right)}{\partial r_i}\right.\]
\[\left.-\frac{\partial \left(\left(\sin(\sum_ir_i)+c_1\right) \left(\left(\sin(\sum_ir_i)+c_1\right)^{-1} - n^{-1}\right)^2 \right)}{\partial r_i}\right)\]
\begin{equation}
-\left(\sin(\sum_ir_i)+c_1\right)g\frac{\partial h_E(\mathbf{r})}{\partial r_j}
\end{equation}
Again, it is taken advantage of the fluid column being split into surface elevation $h_E(x,t)$ and bathymetric fluid depth $h_B(x,t)$

\[-\frac{\partial \sin(\sum_ir_i) }{\partial t} =-\sum_i\left(k_H\frac{\partial \frac{\cos(\sum_ir_i)}{\sin(\sum_ir_i) +c_1} }{\partial r_i}\right.\]
\[\left.+\frac{\partial \left(\frac{\sin(\sum_ir_i)+c_1}{n^2}+ \left(\sin(\sum_ir_i)+c_1\right)^{-1} - \frac{2}{n} \right)}{\partial r_i}\right)\]
\begin{equation}
     -\left(\sin(\sum_ir_i)+c_1\right)g\frac{\partial h(\mathbf{r})-h_B(\mathbf{r})}{\partial r_j}
\end{equation}
\[-\cos\left(\sum_ir_i\right)\]
\[ =\sum_i\left(k_H\frac{\sin(\sum_ir_i)\left(\sin(\sum_ir_i)+c_1\right)+\cos(\sum_ir_i)}{\left(\sin(\sum_ir_i)+c_1\right)^2}\right.\]
\[-\left.\left(\frac{\cos(\sum_ir_i)}{n^2}-\frac{\cos(\sum_ir_i)}{\left(\sin(\sum_ir_i)+c_1\right)^2}\right)\right)\]
\begin{equation}
     +\left(\sin(\sum_ir_i)+c_1\right)g\left(\frac{\partial h_B(\mathbf{r})}{\partial r_j}-\cos\left(\sum_ir_i\right)\right)
\end{equation}
\[\cos\left(\sum_ir_i\right)-\frac{\cos\left(\sum_ir_i\right)}{g\left(\sin(\sum_ir_i)+c_1\right)}\]
\[ -\sum_i\left(\frac{k_H}{g}\left(\frac{\sin(\sum_ir_i)}{\left(\sin(\sum_ir_i)+c_1\right)^2}+\frac{\cos(\sum_ir_i)}{\left(\sin(\sum_ir_i)+c_1\right)^3}\right)\right.\]
\[+\left.\left(\frac{\cos(\sum_ir_i)}{g\left(\sin(\sum_ir_i)+c_1\right)n^2}-\frac{\cos(\sum_ir_i)}{g\left(\sin(\sum_ir_i)+c_1\right)^3}\right)\right)\]
\begin{equation}
     =\frac{\partial h_B(\mathbf{r})}{\partial r_j}
\end{equation}
\[\sin\left(\sum_ir_i\right)-\frac{1}{g}\ln\left(\left|\sin\left(\sum_ir_i\right)+c_1\right|\right)\]
\[ +\sum_i\left(\frac{k_H}{g}\left(a(\mathbf{r})+2^{-1}\left(\sin(\sum_ir_i)+c_1\right)^{-2}\right)\right.\]
\[+\left.\left(\frac{1}{g\,n^2}\ln\left(\left|\sin\left(\sum_ir_i\right)+c_1\right|\right)+(2g)^{-1}\left(\sin(\sum_ir_i)+c_1\right)^{-2}\right)\right)\]
\begin{equation}
     =h_B(\mathbf{r})
\end{equation}
\[\frac{1}{g}\ln\left(\left|\sin\left(\sum_ir_i\right)+c_1\right|\right)\]
\[ -\sum_i\left(\frac{k_H}{g}\left(a(\mathbf{r})+2^{-1}\left(\sin(\sum_ir_i)+c_1\right)^{-2}\right)\right.\]
\[-\left.\left(\frac{1}{g\,n^2}\ln\left(\left|\sin\left(\sum_ir_i\right)+c_1\right|\right)+{2g}^{-1}\left(\sin(\sum_ir_i)+c_1\right)^{-2}\right)\right)\]
\begin{equation}
     =h_E(\mathbf{r})
\end{equation}
and the integral $a(\mathbf{r})$, obtained with Maple, denoted as:
\[a(\mathbf{r})=\frac{2\tan\left(\sum_i\frac{r_i}{2}\right)+2c_1}{\left(c_1^2-1\right)\left(\tan^2\left(\sum_i\frac{r_i}{2}\right)c_1+2\tan\left(\sum_i\frac{r_i}{2}\right)+c_1\right)}\]
\begin{equation}
    +2\arctan\left(\frac{c_1\tan\left(\sum_i\frac{r_i}{2}\right)+1}{\sqrt{c_1^2-1}}\right)\left(c_1^2-1\right)^{-3/2}
    \label{eq:a}
\end{equation}

together with $U_i(\mathbf{r})=\left(\sin(\sum_ir_i)+c_1\right)^{-1}-n^{-1}$ with $i\neq t$, $U_t(\mathbf{r}) =1, and $ $h(\mathbf{r})=c_1+\sin(\sum_ir_i)$. 
\qed
\vspace{4pt}
\newline 
That is, some nonlinear PDE can be transformed into a space such that not the constituents of the PDE itself but the geometries of the space are numerically integrated whereas analytical solutions are maintained for the PDE. This is evident in general relativity where gravitational forces are described as curvature of space. Above solitary solutions permit to validate numerical solvers without having to depict complex domain shapes such as parabolic bowls \cite{bristeau_analytical_2020}. Also, due to the simplicity of the boundary conditions, these can be imposed on the boundaries of arbitrary domain shapes as long as the domain is permitted to adapt along one dimension to $h_B$. For example, any coastal geometry can be used as long as the bathymetry is adapted to $h_B$. This permits the evaluation and assessment of unstructured meshes for the finite volume and finite element method.

\begin{backmatter}


\bmsection{Disclosures} The authors declare no conflicts of interest.

\smallskip

\bmsection{Supplemental document}
See Supplement 1 for supporting content. 

\end{backmatter}

\bigskip

\bibliography{Solitary}

\begin{thebibliography}{10}
\newcommand{\enquote}[1]{``#1''}

\bibitem{Bristeau2018VariousAS}
M.-O. Bristeau, B.~Martino, A.~Mangeney, J.~Sainte-Marie, and F.~Souill{\'e}, \enquote{Various analytical solutions for the incompressible euler and navier-stokes systems with free surface,}  (2018).

\bibitem{sinha_long_2014}
T.~K. Sinha, S.~M.~B. Baruah, and J.~Mathew, \enquote{Long {Wave} {Length} {Soliton} {Solutions} of {Navier} {Stokes} {Equation},} {\protect\JournalTitle{International Journal of Difference Equations}} \textbf{9}, 1--5 (2014).

\bibitem{leoni2020traveling}
G.~Leoni and I.~Tice, \enquote{Traveling wave solutions to the free boundary incompressible navier-stokes equations,}  (2020).

\bibitem{drazin_solitons:_1989}
P.~G. Drazin and R.~S. Johnson, \emph{Solitons: an introduction}, no.~2 in Cambridge texts in applied mathematics (Cambridge University Press, Cambridge [England] ; New York, 1989).

\bibitem{Bredberg2012}
I.~Bredberg, C.~Keeler, V.~Lysov, and A.~Strominger, \enquote{From navier-stokes to einstein,} {\protect\JournalTitle{Journal of High Energy Physics}} \textbf{2012}, 146 (2012).

\bibitem{Padmanabhan:2010rp}
T.~Padmanabhan, \enquote{{Entropy density of spacetime and the Navier-Stokes fluid dynamics of null surfaces},} {\protect\JournalTitle{Phys. Rev. D}} \textbf{83}, 044048 (2011).

\bibitem{Hubeny_2011}
V.~E. Hubeny, \enquote{The fluid/gravity correspondence: a new perspective on the membrane paradigm,} {\protect\JournalTitle{Classical and Quantum Gravity}} \textbf{28}, 114007 (2011).

\bibitem{Rangamani_2009}
M.~Rangamani, \enquote{Gravity and hydrodynamics: lectures on the fluid-gravity correspondence,} {\protect\JournalTitle{Classical and Quantum Gravity}} \textbf{26}, 224003 (2009).

\bibitem{ELING2009496}
C.~Eling, I.~Fouxon, and Y.~Oz, \enquote{The incompressible navier–stokes equations from black hole membrane dynamics,} {\protect\JournalTitle{Physics Letters B}} \textbf{680}, 496 -- 499 (2009).

\bibitem{bristeau_analytical_2020}
M.-O. Bristeau, B.~Di~Martino, A.~Mangeney, J.~Sainte-Marie, and F.~Souillé, \enquote{Some analytical solutions for validation of free surface flow computational codes,} {\protect\JournalTitle{Preprint}}  (2020).

\end{thebibliography}




\end{document}